
\newcount\secno
\newcount\prmno
\def\section#1{\vskip1truecm
               \global\def\currenvir{section}
               \global\advance\secno by1\global\prmno=0
               {\bf \number\secno. {#1}}
               \smallskip}

\def\subsection{\global\def\currenvir{subsection}
                \global\advance\prmno by1
               \smallskip  \ind{ (\number\secno.\number\prmno) }}
\def\subsec{\global\def\currenvir{subsection}
                \global\advance\prmno by1\smallskip
                { (\number\secno.\number\prmno)\ }}

\def\proclaim#1{\global\advance\prmno by 1
                {\bf #1 \the\secno.\the\prmno$.-$ }}

\long\def\th#1 \enonce#2\endth{%
   \medbreak\proclaim{#1}{\it #2}\global\def\currenvir{th}\smallskip}

\def\rem#1{\global\advance\prmno by 1
{\it #1} \the\secno.\the\prmno$.-$ }

\magnification 1250 \pretolerance=500 \tolerance=1000
\brokenpenalty=5000 \mathcode`A="7041 \mathcode`B="7042
\mathcode`C="7043 \mathcode`D="7044 \mathcode`E="7045
\mathcode`F="7046 \mathcode`G="7047 \mathcode`H="7048
\mathcode`I="7049 \mathcode`J="704A \mathcode`K="704B
\mathcode`L="704C \mathcode`M="704D \mathcode`N="704E
\mathcode`O="704F \mathcode`P="7050 \mathcode`Q="7051
\mathcode`R="7052 \mathcode`S="7053 \mathcode`T="7054
\mathcode`U="7055 \mathcode`V="7056 \mathcode`W="7057
\mathcode`X="7058 \mathcode`Y="7059 \mathcode`Z="705A
\def\spacedmath#1{\def\packedmath##1${\bgroup\mathsurround =0pt##1\egroup$}
\mathsurround#1
\everymath={\packedmath}\everydisplay={\mathsurround=0pt}}
 \spacedmath{2pt}

\def\iso{\vbox{\hbox to .8cm{\hfill{$\scriptstyle\sim$}\hfill}
\nointerlineskip\hbox to .8cm{{\hfill$\longrightarrow $\hfill}} }}
\def\sdir_#1^#2{\mathrel{\mathop{\kern0pt\oplus}\limits_{#1}^{#2}}}

\def\pc#1{\tenrm#1\sevenrm}
\def\tx{\kern-1.5pt -}
\def\cqfd{\kern 2truemm\unskip\penalty 500\vrule height 4pt depth 0pt width
4pt\medbreak} 
\def\no{n\up{o}\kern 2pt}
\def\ind{\par\hskip 1truecm\relax}

\font\pal=cmsy7

\def\sp#1{{\cal S}\kern-1pt\raise-1pt\hbox{\pal P}^{}_C(#1)}

\frenchspacing
\input xy
\xyoption{all}
\input amssym.def
\input amssym
\vsize = 25truecm \hsize = 16.1truecm \voffset = -.5truecm
\parindent=0cm
\baselineskip15pt \overfullrule=0pt

\vglue 2.5truecm \font\Bbb=msbm10

\centerline{\bf Class field theory for strictly quasilocal fields}
\smallskip \centerline{\bf with Henselian discrete
valuations\footnote{$^{\ast }$}{\rm Partially supported by Grant
MI-1503/2005 of the Bulgarian Foundation for Scientific Research}}
\bigskip

\centerline{I.D. Chipchakov\footnote{$^{\ast \ast }$}{\rm Mailing
address: I.D. Chipchakov, Institute of Mathematics and
Informatics, Bulgarian Academy of Sciences, Acad. G. Bonchev Str.,
bl. 8, 1113 Sofia, Bulgaria. E-mail: chipchak@math.bas.bg }}

\par
\vskip2.truecm \centerline{{\bf 1. Introduction}}
\par
\medskip
The purpose of this paper is to describe the norm groups of the
fields pointed out in the title. Our starting point is the fact
that a field $K$ is strictly quasilocal, i.e. its finite
extensions are strictly primarily quasilocal (abbreviated,
strictly PQL) fields if and only if these extensions admit
one-dimensional local class field theory (cf. [5, Sect. 3]).
Assuming that $K$ is strictly quasilocal and has a Henselian
discrete valuation $v$, we first show that the norm group
$N(R/K)$ of each finite separable extension $R$ of $K$ is of
index $i(R/K)$ (in the multiplicative group $K ^{\ast }$ of $K$)
dividing the degree $[R\colon K]$. We say that $R$ is a class
field of $N(R/K)$, if $i(R/K) = [R\colon K]$. The present paper
shows that $N(R/K)$ possesses a class field cl$(N(R/K))$ which
is uniquely determined by $N(R/K)$, up-to a $K$-isomorphism. It
proves that cl$(N(R/K))$ includes as a subfield the maximal
abelian extension $R _{\rm ab}$ of $K$ in $R$. Also, we show that
cl$(N(R/K))$ embeds in $R$ as a $K$-subalgebra and is presentable
as a compositum of extensions of $K$ of primary degrees. This
gives rise to a canonical bijection $\omega $ of the set of
isomorphism classes of class fields of $K$ upon the set Nr$(K)$ of
norm groups of finite separable extensions of $K$. Our main
results describe the basic properties of $\omega $ and eventually
enable one to obtain a complete characterization of the elements
of Nr$(K)$ in the set of subgroups of $K ^{\ast }$. They indicate
that $K ^{\ast }$ can be endowed with a structure of a topological
group with respect to which Nr$(K)$ is a system of neighbourhoods
of unity. This topology on $K ^{\ast }$ turns out to be coarser
than the one induced by $v$ unless the residue field $\widehat K$
of $(K, v)$ is finite or of zero characteristic (when they are
equivalent). The present research plays a role in clarifying some
general aspects of one-dimensional abstract local class field
theory, such as the scope of validity of the classical norm
limitation theorem (cf. [11, Ch. 6, Theorem 8]) for strictly
PQL ground fields, and the possibility of reducing the study of
norm groups of quasilocal fields to the special case of finite
abelian extensions (see Remark 3.4, [3] and [6]).
\par
\medskip
The main field-theoretic notions needed for describing the main
results of this paper are the same as those in [6]. Our basic
terminology and notation concerning valuation theory, simple
algebras and Brauer groups, profinite groups, field extensions and
Galois theory are standard (and can be found, for example, in [9;
15; 17] and [21]). As usual, Galois groups are regarded as profinite
with respect to the Krull topology. For convenience of the reader,
we define the notion of a field with (one-dimensional) local class
field theory in Section 2. For each field $E$, $E _{\rm sep}$
denotes a separable closure of $E$, $G _{E} := G(E _{\rm sep}/E)$
is the absolute Galois group of $E$ and $P(E)$ is the set of those
prime numbers $p$ for which $E$ is properly included in its
maximal $p$-extension $E (p)$ in $E _{\rm sep}$. Recall also that
$E$ is said to be PQL, if every cyclic extension $F$ of $E$ is
embeddable as an $E$-subalgebra in each central division
$E$-algebra $D$ of Schur index ind$(D)$ divisible by the degree
$[F\colon E]$. When this occurs, we say that $E$ is a strictly
PQL-field, if the $p$-component Br$(E) _{p}$ of the Brauer group
Br$(E)$ is nontrivial, for every $p \in P(E)$. Let us note that
PQL-fields and quasilocal fields are naturally singled out by the
study of some of the basic types of stable fields with Henselian
valuations (see [1, 5] and the references there). It is also worth
mentioning that they admit satisfactory inner characterizations
which are fairly complete when the considered ground fields belong
to some actively studied and frequently used special classes (see
Section 2 and the observations at the end of [5, Sect. 4], for
more details). The research in this area, however, is primarily
motivated by the fact that strictly PQL-fields admit local class
field theory, and by the validity of the converse in all presently
known cases (cf. [5, Theorem 1 and Sect. 3]). This applies
particularly to the noted place of strictly quasilocal fields in
this theory. As to the choice of our main topic, it is determined
by the fact that the structure of strictly quasilocal fields with
Henselian discrete valuations is known (cf. [4, Sect. 3]) and
sheds light on essential general properties of arbitrary
quasilocal fields. The main results of this paper aim at extending
the traditional basis of one-dimensional abstract local class
field theory (cf. [25] and [2]). They can be stated as follows:
\par
\medskip
{\bf Theorem 1.1.} {\it Let $(K, v)$ be a Henselian discrete
valued strictly quasilocal field with a residue field $\widehat
K$. Then class fields and norm groups of $K$ are related as
follows:}
\par
(i) {\it For each $U \in {\rm Nr}(K)$, there exists a class field
${\rm cl}(U)$ which is uniquely determined, up-to a
$K$-isomorphism; the extension ${\rm cl}(U)/K$ is abelian if and
only if $P(\widehat K)$ contains the prime divisors of the index
of $U$ in $K ^{\ast }$;}
\par
(ii) {\it A class field ${\rm cl}(U)$ of a group $U \in {\rm
Nr}(K)$ embeds as a $K$-subalgebra in a finite extension $R$ of
$K$ in $K _{\rm sep}$ if and only if $N(R/K)$ is included in $U$;
furthermore, if $N(R/K) = U$, then the $K$-isomorphic copy of
${\rm cl}(U)$ in $R$ is unique and includes $R _{\rm ab}$;}
\par
(iii) {\it There exists a set $\{\Phi _{U}\colon \ U \in {\rm
Nr}(K)\}$ of extensions of $K$ in $K _{\rm sep}$, such that $\Phi
_{U}$ is a class field of $U$, $U \in {\rm Nr}(K)$, and for each
$U _{1}, U _{2} \in {\rm Nr}(K)$, $\Phi _{(U _{1} \cap U _{2})}$
equals the compositum $\Phi _{U _{1}}\Phi _{U _{2}}$ and $\Phi
_{(U _{1}U _{2})} = \Phi _{U _{1}} \cap \Phi _{U _{2}}$.}
\par
\vskip0.67truecm {\bf Theorem 1.2.} {\it Assume that $K$, $v$ and
$\widehat K$ satisfy the conditions of Theorem 1.1, ${\rm Op}(K)$ is
the set of open subgroups of $K ^{\ast }$ of finite indices, and
$\Sigma (K)$ is the set of subgroups of $K ^{\ast }$ of finite
indices not divisible by ${\rm char}(\widehat K)$. Then ${\rm
Nr}(K)$, ${\rm Op}(K)$ and $\Sigma (K)$ have the following
properties:}
\par
(i) {\it The intersection of finitely many groups from ${\rm
Nr}(K)$ lies in ${\rm Nr}(K)$; also, if $V$ is a subgroup of $K
^{\ast }$ including a group $U \in {\rm Nr}(K)$, then $V \in {\rm
Nr}(K)$;}
\par
(ii) $K ^{\ast n} \in \Sigma (K)${\it , for each positive integer
$n$ not divisible by {\rm char}$(\widehat K)$;}
\par
(iii) Nr$(K)$ {\it is a subset of ${\rm Op}(K)$ including $\Sigma
(K)$; in order that ${\rm Nr}(K) = {\rm Op}(K)$ it is necessary
and sufficient that ${\rm char}(\widehat K) = 0$ or $\widehat K$
is a finite field.}
\par
\medskip
The paper is organized as follows: Section 2 includes
preliminaries used in the sequel. Section 3 contains the proof of
Theorem 1.1 (i) and a characterization of class fields among
finite separable extensions of a given Henselian discrete valued
strictly quasilocal field. Theorems 1.1 (ii), (iii) and 1.2 are
proved in Section 4.
\par
\vskip1.81truecm \centerline {\bf 2. Preliminaries}
\par
\medskip
Let $E$ be a field, Nr$(E)$ the set of norm groups of finite
extensions of $E$ in $E _{\rm sep}$, and $\Omega (E)$ the set of
finite abelian extensions of $E$ in $E _{\rm sep}$. We say that
$E$ admits (one-dimensional) local class field theory, if the
mapping $\pi $ of $\Omega (E)$ into Nr$(E)$ defined by the rule
$\pi (F) = N(F/E)\colon \ F \in \Omega (E)$, is injective and
satisfies the following two conditions, for each pair $(M _{1}, M
_{2}) \in \Omega (E) \times \Omega (E)$:
\par
The norm group of the compositum $M _{1}M _{2}$ is equal to the
intersection $N(M _{1}/E) \cap $
\par \noindent
$N(M _{2}/E)$ and $N((M _{1} \cap M _{2})/E)$ equals the inner
group product $N(M _{1}/E)N(M _{2}/E)$.
\par
Our approach to the study of fields with such a theory is based on
the following lemma (proved e.g. in [6]).
\par
\medskip
{\bf Lemma 2.1.} {\it Let $E$ be a field and $L/E$ a finite
extension, such that $L = L _{1}L _{2}$, for some extensions $L
_{1}$ and $L _{2}$ of $E$ of relatively prime degrees. Then
$N(L/E) =$
\par \noindent
$N(L _{1}/E) \cap N(L _{2}/E)$, $N(L _{1}/E) = E ^{\ast }
\cap N(L/L _{2})$ and there is a group isomorphism $E ^{\ast
}/N(L/E) \cong (E ^{\ast }/N(L _{1}/E)) \times (E ^{\ast }/N(L
_{2}/E))$.}
\par
\medskip
The following generalization of the norm limitation theorem for
local fields (found in [3]) is used in Section 3 for proving the
existence in the former part of Theorem 1.1 (i).
\par
\medskip
{\bf Proposition 2.2.} {\it Assume that $E$ is a quasilocal field
and $p$ is a prime number, for which the natural Brauer group
homomorphism ${\rm Br}(E) \to {\rm Br}(L)$ maps the $p$-component
${\rm Br}(E) _{p}$ surjectively on ${\rm Br}(L) _{p}$, for every
finite extension $L$ of $E$. Also, let $R/E$ be a finite separable
extension, $R _{\rm ab,p} = R _{\rm ab} \cap E (p)$ and $N _{p}
(R/E)$ the set of those $\alpha \in E ^{\ast }$ for which the
co-set $\alpha N(R/E)$ is a $p$-element in $E ^{\ast }/N(R/E)\}$.
Then $N(R/E) = N _{p} (R/E) \cap N(R _{\rm ab,p}/E)$.}
\par
\medskip
In the rest of this paper, $\overline P$ denotes the set of prime
numbers, and for each field $E$, $P _{0} (E)$ is the subset of
those $p \in \overline P$ for which the polynomial $\sum _{u=0}
^{p-1} X ^{u}$ has a zero in $E$ (i.e. $p = {\rm char}(E)$ or $E$
contains a primitive $p$-th root of unity). Also, we assume that
$P _{1} (E) = \{p ^{\prime } \in (\overline P \setminus P _{0}
(E))\colon \ E ^{\ast } \neq E ^{\ast p'}\}$ and $P _{2} (E) =
\overline P \setminus (P _{0} (E) \cup P _{1} (E))$. Every finite
extension $L$ of a field $K$ with a Henselian discrete valuation
$v$ is considered with its valuation extending $v$, this
prolongation is also denoted by $v$. We write $U(L)$, $\widehat L$
and $e(L/K)$ for the multiplicative group of the valuation ring of
$(L, v)$, the residue field of $(L, v)$ and the ramification index
of $L/K$, respectively. As usual, $U(L) ^{\nu } := \{\lambda ^{\nu
}\colon \ \lambda \in U(L)\}$, for each $\nu \in {\Bbb N}$. The
study of the basic types of PQL-fields with Henselian discrete
valuations relies on the following results (see [2, Sect. 2] and
[4, Sects. 2 and 3]):
\par
\medskip
(2.1) (i) $K$ is PQL if and only if $\widehat K$ is a perfect
PQL-field, and for each $p \in P(\widehat K)$, $\widehat K
(p)/\widehat K$ is a $\hbox{\Bbb Z} _{p}$-extension;
\par
(ii) $K$ is strictly PQL if and only if it is PQL and $P(\widehat
K) = P(K)$;
\par
(iii) $K$ is quasilocal if and only if $\widehat K$ is perfect and
$G _{\widehat K}$ is metabelian of cohomological dimension cd$(G
_{\widehat K}) \le 1$; conversely, every profinite metabelian
group $G$ satisfying the inequality ${\rm cd}(G) \le 1$ is
continuously isomorphic to the absolute Galois group of the
residue field of some Henselian discrete valued quasilocal field
$K(G)$;
\par
(iv) $K$ is strictly quasilocal if and only if the following two
conditions hold: ($\alpha $) $\widehat K$ is perfect with $G
_{\widehat K}$ metabelian of cohomological $p$-dimension cd$_{p}(G
_{K}) = 1$, for each $p \in \overline P$; ($\beta $) $P _{0}
(\widetilde L) \subseteq P(\widetilde L)$, for every finite
extension $\widetilde L$ of $\widehat K$. Conversely, if $G$ is a
metabelian profinite group, such that ${\rm cd} _{p} (G) = 1$, for
all $p \in \overline P$, then $G$ is realizable as an absolute
Galois group of the residue field of a strictly quasilocal field
$K(G)$ with a Henselian discrete valuation.
\par
\medskip
The fulfillment of (2.1) (i) ensures that $K$ is nonreal (cf. [14,
Theorem 3.16]) and the assumption that $K$ is strictly quasilocal
implies the following:
\par
\medskip
(2.2) (i) $P _{0} (K) \setminus \{{\rm char}(\widehat K)\} = P
_{0} (\widehat K) \setminus \{{\rm char}(\widehat K)\}$;
\par
(ii) The group $G _{K}$ is prosolvable (see [2, Proposition 3.1]);
\par
(iii) If $\widetilde L/\widehat K$ is a finite extension, then the
quotient group $\widetilde L ^{\ast }/\widetilde L ^{\ast p ^{\nu
}}$ is cyclic of order $p ^{\nu }$, for every $\nu \in \hbox{\Bbb
N}$ and each $p \in P _{1} (\widehat K)$; this is also true in
case $p \in (P _{0} (\widehat K) \setminus \{{\rm char}(\widehat
K)\})$ and $\nu $ is a positive integer for which $\widehat K$
contains a primitive $p ^{\nu }$-th root of unity;
\par
(iv) Br$(\widetilde L) = \{0\}$ and Br$(L) _{p}$ is isomorphic to
the quasicyclic $p$-group $\hbox{\Bbb Z} (p ^{\infty })$, for every
finite extension $L$ of $K$ and each $p \in P(\widehat K)$ (apply
[21, Ch. II, Proposition 6 (b)] and Scharlau's generalization of
Witt's theorem [19]).
\par
\medskip
It is well-known (cf. [15, Ch. VIII, Sect. 3]) that if
$\varepsilon _{p}$ is a primitive $p$-th root of unity in
$\widehat K _{\rm sep}$, for a given number $p \in (\overline P
\setminus \{{\rm char}(\widehat K)\})$, then the degree $[\widehat
K(\varepsilon _{p})\colon \widehat K]$ divides $p - 1$. Our next result
(proved in [3]) shows that this is the only restriction on the
possible values of the sequence $[\widehat K(\varepsilon
_{p})\colon \widehat K]\colon \ p \in \overline P$. Combined with the
former part of (2.1) (iv), it also describes the behaviour of the
sets $P(\widehat K)$ and $P _{j} (\widehat K)\colon \ j = 0, 1,
2$, when $(K, v)$ runs across the class of Henselian discrete
valued strictly quasilocal fields with char$(\widehat K) = 0$.
\par
\medskip
{\bf Proposition 2.3.} {\it Let $P _{0}$, $P _{1}$, $P _{2}$ and
$P$ be subsets of $\overline P$, such that $P _{0} \cup P _{1}
\cup P _{2}$
\par \noindent
$= \overline P$, $2 \in P _{0}$, $P _{i} \cap P _{j} =
\phi \colon \ 0 \le i < j \le 2$, and $P _{0} \subseteq P
\subseteq (P _{0} \cup P _{2})$. For each $p \in (P _{1} \cup P
_{2})$, let $\gamma _{p}$ be an integer $\ge 2$ dividing $p - 1$
and not divisible by any element of $\overline P \setminus P$.
Assume also that $\gamma _{p} \ge 3$ in case $p \in (P _{2}
\setminus P)$. Then there exists a Henselian discrete valued
strictly quasilocal field $(K, v)$ with the property that $P _{j}
(\widehat K) = P _{j}\colon \ j = 0, 1, 2$, $P(\widehat K) = P$,
and for each $p \in (P _{1} \cup P _{2})$, $\gamma _{p}$ equals
the degree $[K(\varepsilon _{p})\colon K]$, where $\varepsilon
_{p}$ is a primitive $p$-th root of unity in $K _{\rm sep}$.}
\par
\medskip
{\bf Remark 2.4.} It should be pointed out that an abelian torsion
group $T$ is realizable as a Brauer group of a strictly PQL-field
with a Henselian discrete valuation if and only if $T$ is
divisible with a $2$-component $T _{2} \cong \hbox{\Bbb Z} (2
^{\infty })$ (cf. [4, Sect. 4]). This, compared with (2.2) (iv),
shows that strictly PQL-fields form a substantial extension of the
class of strictly quasilocal fields. Such a conclusion can also be
drawn from the study in [7] of the basic types of the considered
fields in the class ${\rm Alg}(E _{0})$ of algebraic extensions of
a global field $E _{0}$. Besides characterizations and a general
existence theorem for such fields, the preprint [7] contains a
description of $N(R/E)$, for a given finite extension $R$ of a
strictly PQL-field $E \in {\rm Alg}(E _{0})$. Similarly to (2.1),
Proposition 2.3 and the results of the present paper, this enables
one not only to find various nontrivial examples of PQL-fields but
also to combine the formal approach to one-dimensional class field
theory (followed in [5], [6, Sect. 3] and [3, Sects. 3-4]) with
efficient constructive methods (see [3, Sect. 5] and [6, Sect.
4]).
\par
\vskip1.81truecm \centerline{\bf 3. Existence and uniqueness of
class fields}
\par
\medskip
The purpose of this Section is to prove Theorem 1.1 (i). Its main
result shows that the validity of the former part of (2.1) (iv)
guarantees the existence of class fields presentable as
compositums of extensions of primary degrees over the ground
fields. It also sheds light on the role of Proposition 2.3 in the
study of norm groups of quasilocal fields.
\par
\medskip
{\bf Theorem 3.1.} {\it Assume that $(K, v)$ is a Henselian
discrete valued strictly quasilocal field and $R$ is a finite
extension of $K$ in $K _{\rm sep}$. Then $R/K$ possesses an
intermediate field $R _{1}$ satisfying the following conditions:}
\par
(i) {\it The sets of prime divisors of $e(R _{1}/K)$, $[\widehat
R _{1}\colon \widehat K]$ and $[\widehat R\colon \widehat R _{1}]$
are included in $P _{1} (\widehat K)$, $\overline P \setminus
P(\widehat K)$ and $P(\widehat K)$, respectively;}
\par
(ii) $N(R/K) = N((R _{\rm ab}R _{1})/K)$ {\it and $K ^{\ast
}/N(R/K)$ is isomorphic to the direct sum $G(R _{\rm ab}/K) \times
(K ^{\ast }/N(R _{1}/K))$;}
\par
(iii) $R _{\rm ab}R _{1}$ {\it is a class field of $N(R/K)$ and
$[(R _{\rm ab}R _{1})\colon K] = [R _{\rm ab}\colon K] \times [R
_{1}\colon K]$.}
\par
\medskip
{\it Proof.} Let $R ^{\prime }$ be the maximal inertial extension
of $K$ in $R$, i.e. the inertial lift of $\widehat R$ in $R$ over
$K$ (cf. [12, Theorems 2.8 and 2.9]). Note first that $R ^{\prime }$
contains as a subfield an extension of $K$ of degree $n _{0}$, for
each $n _{0} \in \hbox{\Bbb N}$ dividing $[R ^{\prime }\colon K]$.
Indeed, by (2.1) (iv), $G _{\widehat K}$ is metabelian with
cd$_{p} (G _{\widehat K}) = 1\colon \ p \in \overline P$, and by
[2, Lemma 1.2], this means that the Sylow pro-$p$-subgroups of $G
_{\widehat K}$ are continuously isomorphic to $\hbox{\Bbb Z}
_{p}$, for each $p \in \overline P$. Therefore, the Sylow
subgroups of the Galois groups of finite Galois extensions of
$\widehat K$ are cyclic. Observing also that if $M$ is the normal
closure of $R ^{\prime }$ in $K _{\rm sep}$ over $K$, then $M/K$
is inertial with $G(M/K) \cong G(\widehat M/\widehat K)$ (cf. [12,
page 135]), one deduces our assertion from Galois theory and the
following lemma.
\par
\medskip
{\bf Lemma 3.2.} {\it Assume that $G$ is a nontrivial finite group
whose Sylow subgroups are cyclic, $H$ is a subgroup of $G$ of
order $n$, and $n _{1}$ is a positive integer divisible by $n$ and
dividing the order $o(G)$ of $G$. Then $G$ possesses a subgroup $H
_{1}$ of order $n _{1}$, such that $H \subseteq H _{1}$.}
\par
\medskip
{\it Proof.} Denote by $p$ the greatest prime divisor of $o(G)$.
Our assumptions show that $G$ is a supersolvable group, whence it
has a normal Sylow $p$-subgroup $G _{p}$ as well as a subgroup $A
_{p}$ isomorphic to $G/G _{p}$ (cf. [13, Ch. 7, Sects. 1 and 2]).
In view of the supersolvability of the subgroups of $G$, this
allows one to prove by induction on $o(G)$ that $G$ has a subgroup
$\widetilde H _{1}$ of order $n _{1}$. Taking finally into account
that $H$ is conjugate in $G$ to a subgroup of $\widetilde H _{1}$
[18] (see also [22, Theorem 18.7]), one completes the proof of the
lemma.
\par
\medskip
Suppose now that $R _{1} ^{\prime }$ is the maximal tamely
ramified extension of $K$ in $R$, $[R _{1} ^{\prime }\colon R
^{\prime }] = n$ and $P$ is the set of prime numbers dividing
$[R\colon K]$. For each $p \in P$, let $f(p)$ and $g(p)$ be the
greatest integers for which $p ^{f(p)} \vert [R\colon K]$ and $p
^{g(p)} \vert [R ^{\prime }\colon K]$. As noted above, Lemma 3.2
indicates that there is an extension $R _{p}$ of $K$ in $R
^{\prime }$ of degree $p ^{g(p)}$, for each $p \in P$. Observing
that $\alpha \in U(R ^{\prime }) ^{n}$, provided that $\alpha \in
R ^{\prime }$ and $v(\alpha - 1) > 0$, one obtains from [16, Ch.
II, Proposition 12] that $R _{1} ^{\prime } = R ^{\prime } (\theta
)$, where $\theta $ is an $n$-th root of $\pi \rho $, for some
$\rho \in U(R ^{\prime })$. Suppose now that $p \in (P _{0}
(\widehat K) \cup P _{1} (\widehat K))$ and $p \neq {\rm
char}(\widehat K)$. Since $p$ does not divide $[R ^{\prime }\colon
R _{p}]$, statement (2.2) (iii) implies the existence of an
element $\rho _{p} \in U(R _{p})$, such that $\rho _{p}\rho ^{-1}$
is a $p ^{(f(p)-g(p))}$-th power in $U(R ^{\prime })$. Therefore,
the binomial $X ^{p ^{(f(p)-g(p))}} - \pi \rho _{p}$ has a zero
$\theta _{p} \in R _{1} ^{\prime }$. Summing up these results, one
proves the following:
\par
\medskip
(3.1) For each $p \in P \cap (P _{0} (\widehat K) \cup P _{1}
(\widehat K))$, $p \neq $ char$(\widehat K)$, there exists an
extension $T _{p}$ of $K$ in $R _{1} ^{\prime }$ of degree $p
^{f(p)}$; moreover, if $p \in P _{0} (\widehat K)$, then the
normal closure of $T _{p}$ in $K _{\rm sep}$ over $K$ is a
$p$-extension.
\par
\medskip
Let $R _{1}$ be the compositum of the fields $R _{p}\colon \ p \in
(P _{2} (\widehat K) \setminus P(\widehat K))$, and $T _{p}\colon \
p \in P _{1} (\widehat K)$. Clearly,  then $R _{1}$ satisfies the
conditions of Theorem 3.1 (i) and embeds as a $K$-subalgebra in the
field $R ^{\prime } (\theta ^{\mu }) := L ^{\prime }$, where $\mu $
is the greatest integer dividing $[R _{1} ^{\prime }\colon R
^{\prime }]$ and not divisible by any element of $\overline P
\setminus P(\widehat K)$. Assuming further that $R _{1} \subseteq L
^{\prime }$, put $L = R ^{\prime }R _{1}$, $T = R _{\rm ab}R _{1}$
and $T ^{\prime } = T(\theta ^{\mu })$. It is easily verified that
$R ^{\prime } \subseteq T$ and $T ^{\prime } = R _{\rm ab}L ^{\prime
}$. As g.c.d.$([R _{\rm ab}\colon K], [R _{1}\colon K]) = 1$, and by
(2.2) (iv), Br$(\widehat K) _{p} \cong \hbox{\Bbb Z} (p ^{\infty
})$, for all $p \in P(\widehat K)$, Lemma 2.1 and [5, Theorem 3.1]
indicate that $K ^{\ast }/N(T/K) \cong G(R _{\rm ab}/K) \times (K
^{\ast }/N(R _{1}/K))$. One also sees that $e(T/K) = e(R _{\rm
ab}/K)e(R _{1}/K)$, $T = R _{\rm ab}L$ and the extensions $L
^{\prime }/L$ and $T ^{\prime }/T$ are tamely totally ramified of
degree $\prod p ^{f(p)-g(p)}$, where $p$ runs through the set $P
\cap (P _{2} (\widehat K) \setminus P(\widehat K))$. In order to
complete the proof of Theorem 3.1, it remains to be shown that
$N(T/K) = N(R/K)$ and $K ^{\ast }/N(R _{1}/K)$ is a group of order
$[R _{1}\colon K]$. Our argument is based on the following two
statements:
\par
\medskip
(3.2) (i) The natural homomorphism of Br$(K) _{p}$ into Br$(Y)
_{p}$ is surjective, for each $p \in P(\widehat K)$ and every
finite extension $Y$ of $K$ in $K _{\rm sep}$;
\par
(ii) $N(T ^{\prime }/K) = N(T/K)$ and $r\pi ^{[R'\colon K]} \in
N(R/K)$, for some $r \in U(K)$.
\par
\medskip
Statement (3.2) (i) is implied by the final assertion of (2.2) (iv)
and the well-known result (cf. [17, Sects. 13.4 and 14.4]) that the
relative Brauer group Br$(Y/E)$ is of exponent dividing $[Y\colon
E]$. The rest of the proof of (3.2) relies on the assumption that
$R ^{\prime }$ is the maximal inertial extension of $K$ in $R$. In
particular, $R$ is totally ramified over $R ^{\prime }$, which
means that $U(R ^{\prime })$ contains an element $\rho $, such
that $\rho \pi \in N(R/R ^{\prime })$. Therefore, the latter part
of (3.2) (ii) applies to the element $r = N _{K} ^{R'} (\rho )$.
In view of (2.1) (iii) and Galois cohomology (cf. [21, Ch. II,
Proposition 6 (b)]), we have $N(\widehat R/\widehat K) = \widehat
K ^{\ast }$, so it follows from the Henselian property of $v$ that
$N(R ^{\prime }/K) = U(K)\langle r\pi ^{[R'\colon K]}\rangle =
U(K)\langle \pi ^{[R'\colon K]}\rangle $. This implies that $N(R
^{\prime }/K)$ is a subgroup of $K ^{\ast }$ of index $[R ^{\prime
}\colon K]$. These observations, combined with the fact that $R
^{\prime } \subseteq T \subseteq T ^{\prime } \subseteq R _{1}
^{\prime }$ and the fields $T, T ^{\prime }$ are tamely and
totally ramified over $R ^{\prime }$, show that $N(T/K) = U(K)
^{e(T/K)}\langle r\pi ^{[R'\colon K]}\rangle $ and $N(T ^{\prime
}/K) = U(K) ^{e(T'/K)}\langle r\pi ^{[R'\colon K]}\rangle $. As
proved above, $[T ^{\prime }\colon T]$ is not divisible by any $p
\in (P (\widehat K) \cup P _{1} (\widehat K))$, whereas
$e(T/K) = e(R _{\rm ab}/K)e(R _{1}/K)$, so it turns out that
g.c.d.$([T ^{\prime }\colon T], e(T/K)) = 1$, $U(K) ^{e(T/K)} =
U(K) ^{e(T'/K)}$ and $N(T/K) = N(T ^{\prime }/K)$. Arguing in a
similar manner, one obtains that $N(R _{1}/K) = U(K) ^{e(R
_{1}/K)}\langle r _{0}\pi ^{[R _{0}\colon K]}\rangle $, where $R
_{0} = R ^{\prime } \cap R _{1}$ and $r _{0}$ is the norm over $K$
of a suitably chosen element of $U(R _{0})$. Since prime divisors
of $e(R _{1}/K)$ lie in $P _{1} (\widehat K)$, the former part of
(2.2) (iii) and the Henselian property of $(K, v)$ imply that
$U(K)/U(K) ^{e(R _{1}/K)}$ is a cyclic group of order $e(R
_{1}/K)$. The obtained results enable one to establish the
required properties of $K ^{\ast }/N(R _{1}/K)$ as a consequence
of the well-known equality $[R _{1}\colon K] = [\widehat R
_{1}\colon \widehat K]e(R _{1}/K)$ (following from Ostrowski's
theorem and from [23, Propositions 2.2 and 3.1]).
\par
It remains to be seen that $N(R/K) = N(T/K)$. The inclusion $N(R/K)
\subseteq N(T/K)$ is evident, so we prove the converse one. Consider
an arbitrary element $\beta $ of $U(K) \cap N(T/K)$, put $[R\colon T
^{\prime }] = m$, and for each $p \in (P \cap P(\widehat K))$, let
$R _{{\rm ab},p} =$
\par \noindent
$R _{\rm ab} \cap K (p)$ and $\rho _{p}$ be the greatest integer
dividing $[R\colon K]$ and not divisible by $p$. It follows from the
inclusion $N(T/K) \subseteq N(R _{{\rm ab},p}/K)$, statement (3.2)
(i) and Proposition 2.2 that $\beta ^{\rho _{p}} \in N(R/K)$, for
each $p \in (P \cap P (\widehat K))$. At the same time, the equality
$N(T/K) = N(T ^{\prime }/K)$ implies that $\beta ^{m} \in N(R/K)$.
Observing now that the prime divisors of $m$ lie in $P(\widehat K)$,
one obtains that g.c.d.$\{m, \rho _{p}\colon \ p \in (P \cap
P(\widehat K))\} = 1$, and therefore, $\beta \in N(R/K)$. Since
$N(T/K) = U(K) ^{e(T/K)}\langle r\pi ^{[R'\colon K]}\rangle $ and
$r\pi ^{[R'\colon K]} \in $
\par \noindent
$N(R/K)$, this means that $N(T/K) \subseteq N(R/K)$, so Theorem 3.1
is proved.
\par
\medskip
The former part of Theorem 1.1 (i) is contained in Theorem 3.1
and the following lemma.
\par
\medskip
{\bf Lemma 3.3.} {\it Let $(K, v)$ be a Henselian discrete
valued strictly quasilocal field, and let $L _{1}$ and $L _{2}$ be
finite extensions of $K$ in $K _{\rm sep}$. Assume also that $L
_{1}$ and $L _{2}$ are class fields of one and the same group $N
\in {\rm Nr}(K)$. Then $L _{1}$ and $L _{2}$ are isomorphic over
$K$.}
\par
\medskip
{\it Proof.} It follows from Theorem 3.1, statement (3.1) and the
observations preceding Lemma 3.2 that one may consider only the
special case in which $[L _{i}\colon K]$ is a $p$-primary number,
for some $p \in (\overline P \setminus P(\widehat K))$. Suppose
first that $p \in (P _{2} (\widehat K) \setminus P(\widehat K))$
and denote by $M$ the minimal normal extension of $K$ in $K _{\rm
sep}$ including $L _{1}$ and $L _{2}$. Then $L _{1}$ and $L _{2}$
are inertial over $K$, whence $M$ has the same property. In view
of (2.1) (iii) and [2, Lemma 1.2], this implies that the Sylow
subgroups of $G(M/K)$ are cyclic, so it follows from [22, Theorem
18.7] that $G(M/L _{1})$ and $G(M/L _{2})$ are conjugate in
$G(M/K)$. Hence, by Galois theory, there exists a $K$-isomorphism
$L _{1} \cong L _{2}$. Assume now that $p \in P _{1} (\widehat
K)$, fix a primitive $p$-th root of unity $\varepsilon \in K _{\rm
sep}$, and put $L ^{\prime } = L(\varepsilon )$, for each finite
extension $L$ of $K$ in $K _{\rm sep}$. It is clear from (2.1)
(iii), the condition on $p$ and the general properties of
cyclotomic extensions (cf. [24, Lemma 1]) that $K ^{\prime }$
contains a primitive $p ^{m}$-th root of unity, for every $m \in
\hbox{\Bbb N}$. Since $[K ^{\prime }\colon K]$ divides $p - 1$
(cf. [15, Ch. VIII, Sect. 3]), our argument also shows that $K
^{\prime \ast } = K ^{\ast }K ^{\prime \ast p ^{m}}$ and $[L _{i}
^{\prime }\colon K ^{\prime }] = [L _{i}\colon K]$. These
observations, combined with Lemma 2.1, imply that the natural
embedding of $K ^{\ast }$ into $K ^{\prime \ast }$ induces an
isomorphism of $K ^{\ast }/N$ on $K ^{\prime \ast }/N(L _{i}
^{\prime }/K ^{\prime })$. As $L _{1}$ and $L _{2}$ are class
fields of $N$, the obtained result leads to the conclusion that
$N(L _{1} ^{\prime }/K ^{\prime }) = N(L _{2} ^{\prime }/K
^{\prime })$. At the same time, (3.1) indicates that the normal
closure of $L _{i} ^{\prime }$ in $K _{\rm sep}$ over $K ^{\prime
}$ is a $p$-extension, so it follows from [6, Theorem 1.1] and the
established properties of $N(L _{i} ^{\prime }/K ^{\prime })$ that
$L _{i} ^{\prime }/K ^{\prime }$ is abelian and $L _{i} ^{\prime
}/K$ is normal. As $K ^{\prime }$ admits local class field theory,
one also sees that $L _{1} ^{\prime } = L _{2} ^{\prime }$. Note
finally that $G(L _{1} ^{\prime }/K)$ is solvable and $G(L _{1}
^{\prime }/L _{1})$ and $G(L _{1} ^{\prime }/L _{2})$ are
subgroups of $G(L _{1} ^{\prime }/K)$ of order $[K ^{\prime
}\colon K]$ and index $[L _{1}\colon K]$. Hence, by P. Hall's
theorem (cf. [13, Ch. 7]), they are conjugate in $G(L _{1}
^{\prime }/K)$, i.e. $L _{1}$ and $L _{2}$ are isomorphic over
$K$, so Lemma 3.3 is proved.
\par
\medskip
{\it Proof of the latter part of Theorem 1.1 (i).} It is clear
from (2.2) (i) and Galois theory that if $A/K$ is a finite
abelian extension, then $[A\colon K]$ is not divisible by any $p
\in (\overline P \setminus P(\widehat K))$. Conversely, let $A$ be
a finite extension of $K$ in $K _{\rm sep}$, such that the prime
divisors of $[A\colon K]$ are contained in $P(\widehat K)$. Then it
follows from (3.2) and Proposition 2.2 that $N(A/E) = N(A _{\rm
ab}/E)$. Thus the latter part of Theorem 1.1 (i) reduces to a
consequence of the former one.
\par
\medskip
{\bf Remark 3.4.} Theorem 3.1 indicates that if $(K, v)$ is a
Henselian discrete valued strictly quasilocal field, $R/K$ is a
finite separable extension and cl$(R/K)$ is the class field of
$N(R/K)$ in $R$, then $R$ is totally ramified over cl$(R/K)$ and
$[R\colon {\rm cl}(R/K)]$ is not divisible by any $p \in P _{1}
(\widehat K)$. In addition, it is easily obtained from (2.1)
(iii), [2, Lemma 1.2] and Galois theory that $\widehat R$ is
abelian over $\widehat K$ if and only if $P(\widehat K)$ contains
all prime divisors of $[\widehat R\colon \widehat K]$. These
results enable one to come to the following conclusions:
\par
(i) $R = {\rm cl}(R/K)$ if and only if the sets of prime divisors
of $[R\colon R _{\rm ab}]$ and $e(R/K)$ are included in $\overline P
\setminus P(\widehat K)$ and $P _{1} (\widehat K) \cup P(\widehat
K)$, respectively;
\par
(ii) $N(R/K) = N(R _{\rm ab}/K)$ if and only if $\widehat
R/\widehat K$ is abelian and $[R\colon K]$ is not divisible by any
$p \in P _{1} (\widehat K)$.
\par
\vskip1.81truecm \centerline{\bf 4. Proofs of Theorems 1.1 (ii),
(iii) and 1.2}
\par
\medskip
Our objective in this Section is to complete the proof of
Theorems 1.1 and 1.2. In what follows, we assume that $K$ is a
strictly quasilocal field with a Henselian discrete valuation
$v$. Note first that the former part of Theorem 1.2 (i) is
implied by Theorem 1.1 (iii). The presentation of the rest of our
argument is divided into three main parts.
\par
\medskip
{\it Proof of Theorem 1.1 (ii) and the latter part of Theorem 1.2
(i).} Let $U$ be a group from Nr$(K)$, $V$ a subgroup
of $K ^{\ast }$ including $U$, $j(u)$ the index of $U$ in $K
^{\ast }$, and $L _{1}$ a class field of $U$ in $K _{\rm sep}$.
Clearly, it is sufficient to prove that $V \in {\rm Nr}(K)$ in
the special case where $j(u)$ is a $p$-primary number. Firstly, if
$p \in P(\widehat K)$, this is implied by Theorem 3.1, Galois
theory and the availability of a local class field theory on $K$.
Secondly, if $p \in (P _{2} (\widehat K) \setminus P(\widehat
K))$, then $L _{1}$ is inertial over $K$ (see Remark 3.4 (i)). As
shown at the beginning of the proof of Theorem 3.1, this indicates
that $K$ has an extension $\Psi _{\tau }$ in $L _{1}$ of degree
$\tau $, for each positive integer $\tau $ dividing $[L _{1}\colon
K]$. Now the assertion that $V \in {\rm Nr}(K)$ follows directly
from the computation of $N(R ^{\prime }/K)$ carried out in the
process of proving (3.2) (ii). Suppose finally that $p \in P _{1}
(\widehat K)$, fix a primitive $p$-th root of unity $\varepsilon
_{p}$ in $K _{\rm sep}$, and put $K ^{\prime } = K(\varepsilon
_{p})$ and $L _{1} ^{\prime } = L _{1} ^{\prime } (\varepsilon
_{p})$. Analyzing the proof of Lemma 3.3, one obtains that there
is a bijection of the set of intermediate fields of $L _{1}
^{\prime }/K ^{\prime }$ on the set of subgroups of $K ^{\prime
\ast }$ including $N(L _{1} ^{\prime }/K ^{\prime })$, and also,
that these intermediate fields are normal over $K$. Thus it turns
out that if $L ^{\prime }$ is the extension of $K ^{\prime }$ in
$L _{1} ^{\prime }$ corresponding to $V.N(L _{1} ^{\prime }/K
^{\prime })$, then $V = N(L/K)$, where $L = L ^{\prime } \cap L
_{1}$. In particular, $V \in {\rm Nr}(K)$, as claimed by the
latter part of Theorem 1.2 (i). In view of Lemma 3.3, this also
proves the former part of Theorem 1.1 (ii).
\par
Assume now that $R$ is a finite extension of $K$ in $K _{\rm sep}$
with $N(R/K) = U$, $P$ is the set of prime numbers dividing
$[R\colon K]$, $R _{1}$ is a subfield of $R$ determined as required
by Theorem 3.1, and $\Phi $ is a class field of $U$ in $R$. By Lemma
3.3, there is a $K$-isomorphism $\Phi \cong R _{\rm ab}R _{1}$, so
$R _{\rm ab}$ is a subfield of $\Phi $. Furthermore, it follows from
Galois theory and the definition of $R _{1}$ that $\Phi $ is
presentable as a compositum of extensions $\Phi _{p}$ of $K$ of
$p$-primary degrees, with $p$ running across $P$. We show that
$\Phi = R _{\rm ab}R _{1}$. It is clearly sufficient to consider
only the special case in which $R \neq R _{\rm ab}$ and to
establish the equality $R _{\rm ab}\Phi _{p} = R _{\rm ab}\Theta
_{p}$, for an arbitrary $p \in (P \setminus P(\widehat K))$ and a
given $K$-isomorphic copy $\Theta _{p}$ of $\Phi _{p}$ in $R _{\rm
ab}R _{1}$. Suppose first that $p \in P _{1} (\widehat K)$. It
follows from the proof of (3.1) that if $R _{\rm ab}\Phi _{p} \neq
R _{\rm ab}\Theta _{p}$, then the set $R \setminus R _{\rm ab}$
must contain a primitive $p$-th root of unity $\varepsilon _{p}$.
However, since the extension $K(\varepsilon _{p})/K$ is abelian
(cf. [15, Ch. VIII, Sect. 3]), this is impossible, so we have $R
_{\rm ab}\Phi _{p} = R _{\rm ab}\Theta _{p}$. Let now $p \in (P
_{2} (\widehat K) \setminus P(\widehat K))$ and $(R _{\rm
ab}\Theta _{p}) \cap (R _{\rm ab}\Phi _{p}) = V _{p}$. In this
case, $\Phi _{p}$ and $\Theta _{p}$ are inertial over $K$, whence
$R _{\rm ab}\Phi _{p}$ and $R _{\rm ab}\Theta _{p}$ are inertial
over $V _{p}$. In view of (2.1) (iv) and Lemma 3.2, this means
that if $R _{\rm ab}\Phi _{p} \neq R _{\rm ab}\Theta _{p}$, then
$R$ possesses distinct subfields $W _{1}$ and $W _{2}$, which are
$V _{p}$-isomorphic extensions of $V _{p}$ of degree $p$. Denote
by $W _{3}$ the minimal normal extension of $V _{p}$ in $K _{\rm
sep}$ including $W _{1}$ and $W _{2}$. Clearly, $W _{3}/V _{p}$ is
a nonabelian
Galois extension, and because of the prosolvability of $G _{K}$
(and the inclusion $(W _{1} \cup W _{2}) \subseteq R$), $W _{3}$
is a subfield of $R$ of degree $pm$ over $V _{p}$, for some
integer $m$ dividing $p - 1$. In addition, it is easily seen that
$W _{3}$ contains as a subfield a cyclic extension of $V _{p}$ of
degree $m$. Observing also that $W _{1}$, $W _{2}$ and $W _{3}$
are inertial over $V _{p}$ (see [12, page 135]), one deduces from
(2.1) (iii) and Galois theory that $m$ is divisible by at least
one number $\mu \in P(\widehat K)$. Thus the hypothesis that $R
_{\rm ab}\Phi _{p} \neq R _{\rm ab}\Theta _{p}$ leads to the
conclusion that $V _{p}$ admits an inertial cyclic extension $Y
_{p}$ in $R$ of degree $\mu $. Since $\mu \in P(\widehat K)$,
statement (2.1) (i) (applied to $K$ and $V _{p}$) implies the
existence of an inertial cyclic extension $Y _{p} ^{\prime }$ of
$K$ with the property that $Y _{p} ^{\prime }V _{p} = Y _{p}$. The
obtained result, however, contradicts the inclusions $Y _{p}
^{\prime } \subseteq R _{\rm ab} \subseteq V _{p}$, so the
equality $R _{\rm ab}\Phi _{p} = R _{\rm ab}\Theta _{p}$ and the
latter part of Theorem 1.1 (ii) are proved.
\par
\medskip
{\it Proof of Theorem 1.1 (iii).} By (2.2) (ii), $G _{K}$ is
prosolvable, whence it possesses a closed Hall pro-$\Pi $-subgroup
$H _{\Pi }$ (uniquely determined, up-to conjugacy in $G _{K}$),
for each subset $\Pi $ of $\overline P$. Denote by $K _{\rm
ur}$ the compositum of inertial finite extensions of $K$ in $K _{\rm
sep}$, and for each $p \in \overline P$, assume that $\Lambda
_{p}$ is the extension of $K$ in $K _{\rm sep}$ corresponding by
Galois theory to a given Hall pro-$(\overline P \setminus
\{p\})$-subgroup of $G _{K}$, $\Omega _{p} = K _{\rm ur} \cap
\Lambda _{p}$, $K _{\rm ab} (p)$ is the maximal abelian extension
of $K$ in $K (p)$, and $N _{p}$ is the set of all groups $X _{p}
\in {\rm Nr}(K)$ of $p$-primary indices in $K ^{\ast }$. Let $I
_{p}$ be the set of finite extensions of $K$ in $\Lambda _{p}$,
$\Omega _{p}$ or $K _{\rm ab} (p)$, depending on whether $p \in P
_{1} (\widehat K)$, $P _{2} (\widehat K) \setminus P(\widehat K)$
or $P(\widehat K)$, respectively. Returning to the proof of Lemma
3.3 and taking into account that $K$ admits local class field
theory (as well as the conjugacy of the Hall pro-$(\overline P
\setminus \{p\})$-subgroups of $G _{K}$), one obtains the
following result:
\par
\medskip
(4.1) The mapping of $I _{p}$ into $N _{p}$, defined by the rule
$\Delta _{p} \to N(\Delta _{p}/K)\colon $
\par \noindent
$\Delta _{p} \in I _{p}$, is bijective, for each $p \in \overline
P$. It transforms field compositums into group intersections and
field intersections into inner group products. Also, every finite
extension of $K$ in $K _{\rm sep}$ of $p$-primary degree is
$K$-isomorphic to a field from $I _{p}$.
\par
\medskip
Suppose now that $U \in {\rm Nr}(K)$, $U \neq K ^{\ast }$, and $P
_{U}$ is the set of prime divisors of the index of $U$ in $K ^{\ast
}$. Then there exists a unique set $\{U _{p}\colon \ p \in P _{U}\}$
of subgroups of $K ^{\ast }$, such that $\cap _{p \in P _{U}} U _{p}
= U$ and each $U _{p}$ is of $p$-primary index in $K ^{\ast }$.
Hence, by the latter part of Theorem 1.2 (i), $U _{p} \in N _{p}$,
and by (4.1), there is a unique field $\Phi _{p} (U) \in I _{p}$
with $N(\Phi _{p} (U)/K) = U _{p}$, for each $p \in P _{U}$.
Denote by $\Phi _{U}$ the compositums of the fields $\Phi _{p}
(U)\colon \ p \in P _{U}$. Applying (4.1) and Lemma 2.1, one
obtains that the set $\{\Phi _{U}\colon \ U \in {\rm Nr}(K)\}$,
where $\Phi _{K ^{\ast }} = K ^{\ast }$, has the properties
required by Theorem 1.1 (iii).
\par
\medskip
{\it Proof of Theorem 1.2 (ii) and (iii).} Let $n$ be a positive
integer not divisible by char$(\widehat K)$, $n _{0}$ the greatest
divisor of $n$ for which $\widehat K$ contains a primitive $n
_{0}$-th root of unity, $n _{1}$ the greatest divisor of $n$ not
divisible by any $p \in (\overline P \setminus P _{1} (\widehat
K))$, and $C _{t}$ a cyclic group of order $t$, for any $t \in
\hbox{\Bbb N}$. It is easily deduced from (2.1) (iv) and (2.2)
(iii) that $\widehat K ^{\ast n} = \widehat K ^{\ast n'}$ and $K
^{\ast }/K ^{\ast n}$ is isomorphic to the direct product $C _{n}
\times C _{n'}$, where $n ^{\prime } = n _{0}n _{1}$. In
particular, it becomes clear that $K ^{\ast n} \in \Sigma (K)$, as
claimed by Theorem 1.2 (ii). One also sees that $K ^{\ast n} =
N((R _{n}T _{n'})/K)$ and $R _{n}T _{n'}$ is a class field of $K
^{\ast n}$, provided that $R _{n}$ and $T _{n'}$ are extensions of
$K$ in $K _{\rm sep}$, such that $R _{n}$ is inertial, $T _{n'}$
is totally ramified, $[R _{n}\colon K] = n$ and $[T _{n'}\colon K]
= n ^{\prime }$. This, combined with Theorem 1.2 (i), proves that
$\Sigma (K) \subseteq {\rm Nr}(K)$. As to the inclusion Nr$(K)
\subseteq {\rm Op}(K)$, it can be viewed as a consequence of
Theorem 3.1, since it is well-known that the groups from Nr$(K)$
are open in $K ^{\ast }$. Note also that the classical existence
theorem yields Nr$(K) = {\rm Op} (K)$ in case $\widehat K$ is a
finite field (cf. [11, Ch. 6] or [9, Ch. IV]), so it remains to be
seen that Nr$(K) \neq {\rm Op}(K)$ whenever $\widehat K$ is
infinite of characteristic $q > 0$ and cardinality $\kappa $.
Arguing as in the proof of Propositions 3 and 4 of [25, Part IV]
(or of [9, Ch. V, (3.6)]), one obtains this result by showing that
Nr$(K)$ and the set of finite abelian extensions of $K$ in $K (q)$
are of cardinality $\kappa $ whereas Op$(K)$ is of cardinality $2
^{\kappa }$. Thus Theorems 1.1 and 1.2 are proved.
\par
\medskip
{\bf Remark 4.1.} It is easily obtained from (2.2) (iii) that if $U$
is a subgroup of $K ^{\ast }$ of finite index $n$ not divisible by
char$(\widehat K)$, then $K ^{\ast }/U$ is isomorphic to $C _{e}
\times C _{n/e}$ and $n$ divides $e ^{2}$, where $e$ is the
exponent of $K ^{\ast }/U$.
\par
\medskip
Note finally that Theorem 1.2 fully characterizes the elements of
Nr$(K)$ in the set of subgroups of $K ^{\ast }$, provided that
char$(\widehat K) = 0$. If char$(\widehat K) = q > 0$ and $L$ is a
finite extension of $K$ in $K _{\rm sep}$, then Theorem 3.1 and
Lemma 2.1 imply that $N(L/K) = N(L _{\rm ab,q}/K) \cap N(L
_{0}/K)$, for some extension $L _{0}$ of $K$ in $L$ of degree not
divisible by $q$. This, combined with Hazewinkel's existence
theorem [10] concerning totally ramified abelian $q$-extensions of
$K$ (see also [8, 3.5 and 3.7] and [20]), allows one to obtain a
satisfactory inner characterization of the groups from Nr$(K)$.
\par
\medskip
Acknowledgements. The results of Section 3 were obtained during my
visit to Tokai University, Hiratsuka, Japan (9. 2002/ 3. 2003). I
gratefully acknowledge the stimulating atmosphere and the
excellent conditions at the University (including the ensured
efficient fully supported medical care after my involvement in a
travel accident). It is a pleasure for me to thank my
host-professor M. Tanaka, Mrs. M. Suzuki, Mrs. A. Uchida and the
colleagues at the Department of Mathematics for their kind
hospitality. My thanks are also due to Grant MI-1503/2005 of the
Bulgarian Foundation for Scientific Research for the partial
support and to the referee for his or her suggestions.
\vskip1cm \centerline{ REFERENCES} \vglue15pt\baselineskip12.8pt
\def\num#1{\smallskip\item{\hbox to\parindent{\enskip [#1]\hfill}}}
\parindent=1.38cm
\par
\num{1} {\pc CHIPCHAKOV}, I.D.: {\sl Henselian valued stable fields.} J.
Algebra {\bf 208}, 344-369 (1998).
\par
\num{2} {\pc CHIPCHAKOV}, I.D.: {\sl Henselian valued quasilocal fields
with totally indivisible value groups.} Comm. Algebra {\bf 27},
3093-3108 (1999).
\par
\num{3} {\pc CHIPCHAKOV}, I.D.: {\sl On the scope of validity of the norm
limitation theorem for quasilocal fields.} Preprint (available at
www.arXiv.org/RA.math/0511587).
\par
\num{4} {\pc CHIPCHAKOV}, I.D.: {\sl Henselian discrete valued fields
admitting one-dimensional local class field theory.} In:
Proceedings of AGAAP-Conference (V. Brinzanescu, V. Drensky and
P. Pragacz, Eds.), 23.9-02.10. 2003, Borovets, Bulgaria, Serdica
Math. J. {\bf 30}, 363-394 (2004).
\par
\num{5} {\pc CHIPCHAKOV}, I.D.: {\sl One-dimensional abstract local class
field theory.} Preprint (available at www.arXiv.org/
RA.math/0506515).
\par
\num{6} {\pc CHIPCHAKOV}, I.D.: {\sl On nilpotent Galois groups and
the scope of the norm limitation theorem in one-dimensional abstract
local class field theory.} Preprint (to appear in: Proceedings of
ICTAMI 05, 15.9-18.9 2005, Albac, Romania, Acta Universitatis
Apulensis {\bf 10}).
\par
\num{7} {\pc CHIPCHAKOV}, Algebraic extensions of global fields
admitting one-dimensional local class field theory. Preprint
(available at www.arXiv.org/
\par
NT.math/0504021).
\par
\num{8} {\pc FESENKO}, I.B.: {\sl Local class field theory.} Perfect residue
field case. Izv. Ross. Akad. Nauk, Ser. Mat. {\bf 57}, No 4, 72-91
(1993) (Russian: Engl. transl. in Russ. Acad. Sci., Izv., Math.
{\bf 43}, No 4, 72-91 (1993).
\par
\num{9} {\pc FESENKO}, I.B., {\pc VOSTOKOV}, S.V.: {\sl Local Fields and Their
Extensions.} Transl. Math. Monographs, 121, Am. Math. Soc.,
Providence, RI, 2002.
\par
\num{10} {\pc HAZEWINKEL}, M.: {\sl Corps de classes local.} Appendix to {\pc DEMAZURE}, M., {\pc GABRIEL}, P.: {\sl Groupes Algebriques,} North-Holland,
Amsterdam, 1970.
\par
\num{11} {\pc IWASAWA}, K.: {\sl Local Class Field Theory.} Iwanami Shoten,
Japan, 1980 (Japanese: Russian transl. in Mir, Moscow, 1983).
\par
\num{12} {\pc JACOB}, B., {\pc WADSWORTH}, A.: {\sl Division  algebras  over
Henselian fields.} J. Algebra {\bf 128}, 126-179 (1990).
\par
\num{13} {\pc KARGAPOLOV}, M.I., {\pc MERZLYAKOV}, Yu.I.: {\sl Fundamentals of the
Theory of Groups, 2nd ed.,} Nauka, Moscow, 1977 (in Russian: Engl.
transl. in Graduate Texts in Mathematics {\bf 62},
Springer-Verlag, New York-Heidelberg-Berlin, 1979).
\par
\num{14} {\pc LAM}, T.Y.: {\sl Orderings, valuations and quadratic
forms.} Reg. Conf. Ser. Math. {\bf 52}, 1983.
\par
\num{15} {\pc LANG}, S.: {\sl Algebra.} Addison-Wesley Publ. Comp., Mass.,
1965.
\par
\num{16} {\pc LANG}, S.: {\sl Algebraic Number Theory.} Addison-Wesley Publ.
Comp., Mass., 1970.
\par
\num{17} {\pc PIERCE}, R.: {\sl Associative Algebras.} Springer-Verlag, New
York, 1982.
\par
\num{18} {\pc RUSAKOV}, S.A.: {\sl Analogues to Sylow's theorems on
existence and embeddability of subgroups.} Sibirsk. Mat. Zh. {\bf
4}, No 5, 325-342 (1963) (in Russian).
\par
\num{19} {\pc SCHARLAU}, W.: {\sl $\ddot U$ber die Brauer-Gruppe eines
Hensel-K$\ddot o$rpers.} Abh. Math. Semin. Univ. Hamb. {\bf 33},
243-249 (1969).
\par
\num{20} {\pc SEKIGUCHI}, K.: {\sl Class field theory of $p$-extensions over a
formal power series field with a $p$-quasifinite coefficient field.}
Tokyo J. Math. {\bf 6}, 167-190 (1983).
\par
\num{21} {\pc SERRE}, J.-P.: {\sl Cohomologie Galoisienne.} Lecture Notes in
Math. {\bf 5}, Springer-Verlag, Berlin-Heidelberg-New York, 1965.
\par
\num{22} {\pc SHEMETKOV}, L.A.: {\sl Formations of Finite Groups.} Nauka,
Moscow, 1978 (in Russian).
\par
\num{23} {\pc TOMCHIN}, T.Y., {\pc YANCHEVSKIJ}, V.I.: {\sl On defects of valued
division algebras.} Algebra I Analiz {\bf 3}, 147-164 (1991) (in
Russian: Engl. Transl. in St. Petersburg Math. J. {\bf 3}, 631-646
(1992)).
\par
\num{24} {\pc WARE}, R.: {\sl Galois groups of maximal $p$-extensions.} Trans. Am.
Math. Soc. {\bf 333}, 721-728 (1992).
\par
\num{25} {\pc WHAPLES}, G.: {\sl Generalized local class field theory. II.
Existence Theorem. III. Second form of the existence theorem.
Structure of analytic groups. IV. Cardinalities.} Duke Math. J. {\bf
21}, 247-255, 575-581, 583-586 (1954).
\par
\end